\documentclass[12pt]{article}
\usepackage{amsmath,amssymb}
\usepackage{titletoc}
\usepackage{latexsym}
\usepackage{amsmath}
\usepackage{amssymb}
\usepackage{multicol}
\usepackage{graphics}
\usepackage{graphicx}
\usepackage{subfigure}
\usepackage{indentfirst}
\usepackage{epsfig}
\usepackage{mathrsfs}
\newcommand{\p}{\partial}

\newtheorem{thm}{Theorem}
\newtheorem{cor}{Corollary}
\newtheorem{lem}{Lemma}

\title{Uniqueness in determination of the fractional order in the TFDE using one measurement
\author{Yi Zhang \hspace{2mm} Xianzheng Jia \hspace{2mm} Gongsheng Li{\footnote{Corresponding author, Email: ligs@sdut.edu.cn.}}\\
\emph{\small School of Mathematics and Statistics, Shandong University of Technology}\\
\emph{\small Zibo, Shandong 255049, PRC}
}
\date{}
}

\begin{document}

\maketitle
\begin{center}
\begin{minipage}{13cm}
\vskip 0.2cm
{\bf Abstract}: This article deals with an inverse problem of identifying the fractional order in the 1D time fractional diffusion equation (TFDE in short) using the measurement at one space-time point. Based on the expression of the solution to the forward problem, the inverse problem is transformed to a nonlinear algebraic equation. By choosing suitable initial values and the measured point, the nonlinear equation has a unique solution by the monotonicity of the Mittag-Lellfer function. Theoretical testifications are presented to demonstrate the unique solvability of the inverse problem.
\end{minipage}
\end{center}

\section{Introduction}
Many diffusion processes in nature and engineering, such as contaminants transport in the soil, oil flow in porous media,
long distance transport of pollutants in the groundwater, etc., are referred to as anomalous diffusion, where the particle plume spreads slower or faster than predicted by the classical integer-order diffusion model. In recent three decades, fractional diffusion equations have been found to be efficient mathematical models to describe some anomalous diffusion phenomena (see \cite{ada92, ben98, ber00, hat98, main10, xiong06, zhou03}, e.g.). The TFDE model, that is obtained from the classical diffusion equation by replacing the first-order time derivative by a fractional derivative, is given as
$$
\p_t^\alpha u - D \triangle u = 0, (x, t)\in \Omega_T, \eqno{(1.1)}
$$
where $\Omega_T=\Omega\times (0, T)$ with $T>0$, $\Omega\subset {\bf R}^N$ is an open bounded domain with smooth boundary, and $u=u(x,t)$ denotes the state variable at space point $x$ and time $t$, $D>0$ is the diffusion coefficient, $\alpha\in (0, 1)$ is the fractional order, and the fractional derivative $\p_t^\alpha u$ is defined by Caputo's definition:
$$
\p_t^\alpha u=\frac{1}{\Gamma(1-\alpha)}\int_{0}^{t}\frac{\p u(x,s)} {\p s}\frac{ds}{(t-s)^{\alpha}}.
\eqno{(1.2)}
$$
See, e.g., Podlubny \cite{pod99} and Kilbas et al \cite{kil06} for the definitions and properties of the fractional derivatives.\\
\indent There are quite a lot of researches on equation (1.1) from numerics, and some works from theoretical aspects, we refer to \cite{chen12, ding15, eide04, gor99, han02, luch09, luch10, metz00, saka11} for the theoretical analysis on the solution of the forward problem, and also see the monographs \cite{pod99}, \cite{kil06}, \cite{kub20} and \cite{gor20}. On the other hand, we refer to \cite{ali15, chi11, gue21, jia14, jin15, liz16, liuj10, liuy16, mur07, weit16,yama12} for the researches on inverse problems arising from the TFDEs. As we know the fractional order in a fractional model is a key parameter to characterize the heavy-tail subdiffusion with memory effect. However, it is always unknown for real-life problems which leading to inverse fractional order problems. The inverse problems of identifying parameters including fractional orders have been studied during the last decade. We refer to Alimor and Ashurov \cite{ali20}, Ashurov and Umarov \cite{ash201}, Ashurov and Zunnunov \cite{ash202}, Cheng et al \cite{cheng09}, Hatano et al \cite{hat13}, Janno \cite{jan16}, Janno and Kinash \cite{jan18}, Jin and Kian \cite{jin21}, Li et al \cite{li13}, Li et al \cite{liz20}, Li and Yamamoto \cite{liz151}, Sun et al \cite{sunl21}, Tatar and Ulusoy \cite{tatar13}, Yamamoto \cite{yama20,yama21}, Yu et al \cite{yu15}, and so on. \\
\indent It is noted that in the existing work on inverse fractional order problems, most of them were studied by using subdomain measurements or one-point measurements at $t\in (0, T)$, or using subboundary data also at $t\in (0, T)$ for arbiutrary given $T>0$. An interesting problem is can we determine the fractional orders in theory only using limited discretized data? Exactly speaking, there is one order $\alpha\in (0, 1)$ in equation (1.1) which is unknown, can we determine it uniquely only using one measurement?\\
\indent It is difficult to give an answer in theory for the above question, but the situation could be changed if having suitable conditions. Here we are concerned with the inverse problem of determining the fractional order $\alpha$ in Eq. (1.1) using one measurement at one space-time point. By the eigenfunction expansion method, the solution of the forward problem is expressed by the Mittag-Lellfer function, and the inverse problem is transformed to a nonlinear algebraic equation. By choosing suitable initial values and the measured point, the nonlinear equation can be solved uniquely by the monotonicity of the Mittag-Lellfer function. The unique solvability of the inverse problem is testified by theoretical examples.\\
\indent The rest of the paper is organized as follows.\\
\indent In Section 2, some preliminaries on the Mittag-Lellfer function, and solvability on a nonlinear algebraic equation are introduced. The forward problem and its solution, and the inverse problem of identifying the fractional order using one measurement are given in Section 3. In Section 4, a local well-posedness of the inverse problem is obtained by the theory of solution to the nonlinear algebraic equation, and a conditional uniqueness is proved by suitably choosing the initial value and the measured time, and numerical testifications are presented. Conclusions are given in Section 5.

\section{Preliminaries}
We give some basic facts on the Gamma function, the Mittag-Leffler function and their properties, and the theory of solving a nonlinear algebraic equation. We refer to the monographs \cite{gor20} and \cite{liq05} for details.
\subsection{The Gamma function and Mittag-Leffler function}
The Gamma function is an analytic continuation of the factorial function in the entire complex plane, which is defined by
$$
\Gamma(z)=\int_0^\infty x^z e^{-x}\ dx, \eqno{(2.1)}
$$
with $\Re(z)>0$. For the Gamma function there hold the formula $\Gamma(z+1)=z \Gamma(z)$, and the Stirling's approximate formula
$$
\Gamma(z+1)\sim \sqrt{2\pi z}\ e^{-z}\ z^z,\ z\rightarrow +\infty, \eqno{(2.2)}
$$
and the derivative's formula
$$
\frac{{{\Gamma^{'}}(z)}}{{\Gamma(z)}} =  - \gamma - \frac{1}{z} + \sum\limits_{n = 1}^\infty {(\frac{1}{n} - \frac{1}{{n + z}})},
\eqno{(2.3)}
$$
where ${\Gamma^{'}}(z)$ denotes the derivative of the Gamma function, and $\gamma=0.5772157\dots$ is the  Euler constant. On the Gamma function, we have the following assertion.
\begin{lem}
For the parameter $\alpha\in (0, 1)$, there holds
$$
\lim\limits_{j\rightarrow \infty} \frac{\Gamma(\alpha j)}{\Gamma(\alpha j+\alpha)} =0.\eqno{(2.4)}
$$
\end{lem}
{\emph{Proof}} By (2.2) we have for large $z>0$
$$
\begin{array}{lll}
\frac{\Gamma(\alpha j)}{\Gamma(\alpha j+\alpha)} &\sim & \frac{{{{[\sqrt{2\pi} {e^{-\alpha j}}{{(\alpha j)}^{\alpha j}}]} \mathord{\left/
 {\vphantom {{[\sqrt{2\pi} {e^{-\alpha j}}{{(\alpha j)}^{\alpha j}}]} {\sqrt{\alpha j}}}} \right.
 \kern-\nulldelimiterspace} {\sqrt{\alpha j}}}}}{{{{[\sqrt{2\pi} {e^{-\alpha j - \alpha}}{{(\alpha j + \alpha)}^{\alpha j + \alpha}}]} \mathord{\left/
 {\vphantom {{[\sqrt{2\pi} {e^{-\alpha j - \alpha}}{{(\alpha j + \alpha)}^{\alpha j + \alpha}}]} {\sqrt{\alpha j + \alpha}}}}\right.
 \kern-\nulldelimiterspace} {\sqrt{\alpha j + \alpha}}}}}\ (j \rightarrow \infty)\\
 & \sim & {e^\alpha}\  {\left({\frac{{\alpha j}}{{\alpha j + \alpha}}}\right)^{\alpha j}}\  \frac{1}{{{{(\alpha j + \alpha)}^\alpha}}}\  \sqrt{\frac{{\alpha j + \alpha}}{{\alpha j}}}\ (j \rightarrow \infty)\\
 & \sim & {e^\alpha}\  {\left[{{{\left({\frac{1}{{1 + {1 \mathord{\left/
 {\vphantom {1 j}}\right.
 \kern-\nulldelimiterspace} j}}}}\right)}^j}}\right]^\alpha}\  \frac{1}{{{{(\alpha j + \alpha )}^\alpha}}}\  \sqrt{1 + \frac{1}{j}}\ (j \rightarrow \infty).\\
\end{array}
\eqno{(2.5)}
$$
Thanks to $\mathop {\lim }\limits_{j \to \infty } {\left( {\frac{1}{{1 + {1 \mathord{\left/{\vphantom {1 j}} \right.\kern-\nulldelimiterspace} j}}}} \right)^j} = {e^{-1}}$, we get
$$
\begin{array}{lll}
\mathop{\lim}\limits_{j \to \infty} \frac{\Gamma(\alpha j)}{\Gamma(\alpha j+\alpha)}
&=& \mathop{\lim}\limits_{j \to \infty} \left[{{e^\alpha}\  {e^{-\alpha}}\  \frac{1}{{{{(\alpha j + \alpha)}^\alpha}}} \ \sqrt{1 + \frac{1}{j}}}\right]\\
&=& \mathop{\lim}\limits_{j \to \infty} \frac{1}{(\alpha j + \alpha)^\alpha}=0.
\end{array}
\eqno{(2.6)}
$$
\indent Following Mittag-Leffler's classical definition, the one-parametric Mittag-Leffler function is defined by the power series
$$
E_{\alpha}(z) = \sum\limits_{j = 0}^\infty  {\frac{{{z^j}}}{{\Gamma(\alpha j + 1)}}},\ z\in {\cal C},\ \alpha>0.
\eqno{(2.7)}
$$
Obviously there is $E_1(z)=e^z$ as $\alpha=1$, i.e., $E_{\alpha}(z)$ is a generalization of the exponential function $e^z$.\\
\indent On analytical properties of the Mittag-Leffler functions, we have the following assertions.
\begin{lem}
(i) There exists a constant $c>0$ such that $|E_{\alpha}(z)|\leq \frac{c}{1+|z|}$, $\alpha>0$.\\
(ii) For real number $z\in {\bf R}$, the function $E_{\alpha}(z)$ is strictly monotone on $z>0$, and decreasing on $z<0$.
\end{lem}
\begin{lem}
Let $G(\alpha)=E_\alpha(-c t^\alpha)$ for $\alpha\in (0, 1)$ and given $0<c, t<\infty$. There holds
$$
G'(\alpha)=\sum\limits_{j=1}^\infty (-c)^j\ j\ t^{\alpha j} \frac{\ln(t)+\gamma_j}{\Gamma(\alpha j+1)}, \eqno{(2.8)}
$$
where $\gamma_j=\gamma + \frac{1}{\alpha j+1} - \sum\limits_{n = 1}^\infty {(\frac{1}{n} - \frac{1}{{n + \alpha j+1}})}$ according to (2.3).
\end{lem}
{\emph{Proof}}\quad By (2.7), it is easy to deduce that the derivative $G'(\alpha)$ ($0<\alpha<1$) is expressed by (2.8), where $\gamma_j$ is given by (2.3), and we only need to prove the convergence of the series at the right-hand side of (2.8).\\
\indent Denote
$$
y_j=\frac{(-c)^j\ j\ t^{\alpha j}}{\Gamma(\alpha j+1)} \ln(t),\ \bar{y}_j=\frac{(-c)^j\ j\ t^{\alpha j}\gamma_j}{\Gamma(\alpha j+1)},\  j=1, 2, \cdots,\eqno{(2.9)}
$$
and we have $G'(\alpha)=\sum\limits_{j=1}^\infty (y_j+\bar{y}_j)$. For given $c$ and $t$ as parameters, there holds
$$
\left|\frac{y_{j+1}}{y_j}\right|=c\  t^\alpha \frac{(j+1)\  \Gamma(\alpha j+1)}{j\  \Gamma(\alpha j+\alpha+1)}=c\  t^\alpha \frac{\Gamma(\alpha j)}{\Gamma(\alpha j+\alpha)},\eqno{(2.10)}
$$
here we utilize the formula $\Gamma(z+1)=z \Gamma(z)$. Then by Lemma 1 we get
$$
\lim\limits_{j\rightarrow \infty} |\frac{y_{j+1}}{y_j}|=0,\eqno{(2.11)}
$$
by which the series $\sum\limits_{j=1}^\infty y_j$ is convergent. Noting that
$$
\lim\limits_{j\rightarrow \infty} \frac{\gamma_{j+1}}{\gamma_j} = 1,\eqno{(2.12)}
$$
we deduce that $\sum\limits_{j=1}^\infty \bar{y}_j$ is also convergent by the same arguments. So the derivative $G'(\alpha)$ is well-defined for $\alpha\in (0, 1)$.
\subsection{Solution to a nonlinear algebraic equation}
Consider to solve a group of nonlinear algebraic equations
$$
F(x)=y,\eqno{(2.13)}
$$
where $F: D\subset {\bf R}^N\rightarrow {\bf R}^N$ is continuous differential on $D$, and $D$ is an open bounded domain in ${\bf R}^N$. On the local unique solvability of the equation (2.13), we give the following assertion \cite{liq05}.
\begin{lem}
For the equation (2.13), assume that there exists a nonsingular matrix $B\in {\cal L}({\bf R}^N)$, and $x^0\in D$ with a closed sphere $S_0=\bar{S}(x^0, \delta)$ such that
$$
||F(x)-F(y)-B(x-y)||\leq \beta ||x-y||,\ \forall x, y\in S_0,\eqno{(2.14)}
$$
where $0<\beta<\|B^{-1}\|^{-1}$. Then $F$ is a homeomorphsim mapping from $S_0$ to $F(S_0)$.
\end{lem}

\indent We now consider to solve a nonlinear algebraic equation in 1D case
$$
g(x)=y,\eqno{(2.15)}
$$
where $g: D\subset {\bf R}\rightarrow {\bf R}$ is continuous differential on $D$, and $D$ is an open bounded interval in ${\bf R}$. Based on Lemma 4 there holds
\begin{cor}
Suppose that the function $g: D\subset {\bf R}\rightarrow {\bf R}$ is continuous differentiable at $x_0\in int(D)$, and $g'(x_0)\neq 0$, and $V_0=\{x: |x-x_0|\leq \delta\}$ for any $\delta>0$ as a neighborhood of $x_0$. Then $g$ is a local homeomorphism mapping at $x_0$, i.e., the equation $g(x)=y$ has a unique solution for given $y\in g(V_0)$. Moreover, the inverse function $g^{-1}: g(V_0)\rightarrow V_0$ is also continuous.
\end{cor}
{\emph{Proof}}\quad Denote $a=g'(x_0)$ and $a\neq 0$. By the definition of the derivative, we have
$$
\lim\limits_{x\rightarrow x_0} \frac{g(x)-g(x_0)}{x-x_0} =a.\eqno{(2.16)}
$$
Thus $\forall \varepsilon>0$ and $\varepsilon<|a|$, there exists $\delta>0$ such that
$$
|g(x)-g(x_0) - a (x-x_0)|\leq \varepsilon |x-x_0|,\  \forall x, x_0\in V_0.\eqno{(2.17)}
$$
Noting $\varepsilon<|a|$, by Lemma 4 follows that the mapping $g$ is a local homeomorphism from $V_0$ to $g(V_0)$, and the equation $g(x)=y$ has a unique solution for given $y\in g(V_0)$ and the inverse function $g^{-1}$ is continuous.
\section{The forward problem and the inverse problem}
\subsection{The forward problem}
Consider the homogeneous time-fractional diffusion equation in 1D case
$$
\p_t^\alpha u - D u_{xx} = 0, (x, t)\in \Omega_T, \eqno{(3.1)}
$$
where $\Omega_T=\Omega\times (0, T)$ and $\Omega=(0, l)$ for $l>0$. For the equation (3.1), given the homogeneous Dirichelet boundary condition
$$
u(0, t)=u(l,t)=0,\ 0<t\leq T, \eqno{(3.2)}
$$
and the initial distribution
$$
u(x, 0)=f(x),\ x\in \Omega, \eqno{(3.3)}
$$
where $f(x)\in L^2(\Omega)$. A forward problem is composed by the equation (3.1) with the initial boundary value conditions (3.2)-(3.3). By using the eigenfunction expansion method, there exists a unique solution (see \cite{saka11,yama21} for instance):
$$
u(x,t)\in C([0, T], L^2(\Omega))\cap C((0, T], H^2(\Omega)\cap H_0^1(\Omega)),\eqno{(3.4)}
$$
and it can be expressed by the Mittag-Lellfer function given as
$$
u(x, t)=\sum\limits_{n = 1}^\infty f_n\ {E_{\alpha}} (-D {\lambda_n}{t^\alpha})\ {\varphi_n}(x), \eqno{(3.5)}
$$
where $f_n=\frac{2}{l} (f, \varphi_n)$, and $\lambda_n=\frac{n^2 \pi^2}{l^2}$, $\varphi_n(x)=\sin(\frac{n\pi x}{l})$, and $E_{\alpha}(\cdot)$ is given by (2.7).
\subsection{The inverse problem}
Although we get the solution to the forward problem given by (3.5), it can not be put into practice if there are unknown parameters in the model. Actually, the fractional order $\alpha$ in Eq. (3.1) is an important index characterizing the anomalous diffusion with memory effect in time, but it cannot be measured which leads to inverse fractional order problems.\\
\indent Suppose that the fractional order in Eq. (3.1) is unknown, and we are to determine it by limited data in time at one space point. Let $x_0\in \Omega$ be fixed, and we have the measured data at time $t_k$ given as
$$
u(x_0, t_k),\  k=1,2,\cdots, K,\eqno{(3.6)}
$$
where $K\geq 1$. The inverse problem is to identify the order $\alpha\in (0, 1)$ using the data (3.6) based on the forward problem (3.1)-(3.3).\\
\indent Denote $I=(0, 1)$, and $u(\alpha)(x, t)$ as the unique solution to the forward problem for any given $\alpha\in I$. An interesting problem is:\\
\indent{\it Can we determine $\alpha\in I$ uniquely only using one measurement of $u(x_0, t_1)$ at $t_1\in (0, T)$} ?\\
\indent Generally speaking, it is difficult to answer the above question in theory. Nevertheless, we can determine it in some special cases. Denote $d=u(x_0, t_1)$ as the measurement. Then noting the solution's expression (3.5), we get a nonlinear algebraic equation:
$$
F(\alpha)=u(\alpha)(x_0,t_1)=\sum\limits_{n = 1}^\infty f_n\ {E_{\alpha}} (-D {\lambda_{n}}{t_1^\alpha})\ {\varphi_n}(x_0)=d.
\eqno{(3.7)}
$$
As a result the inverse fractional order problem is transformed to solving of the nonlinear equation (3.7).

\section{Unique solvability of the inverse problem}
\subsection{Local well-posedness}
We now prove the inverse problem is locally solvable and has unique solution for $\alpha\in I$ by choosing appropriate measured point and the measured time $t_1>0$. For that we need the existence of the derivative of the function $F(\alpha)$ on $\alpha\in I$, and some limitations for the initial value functions.\\
\indent Without loss of generality, we set $l=\pi$. It is known that the eigenfunction system of the forward problem is $\{n^2, \sin(nx)\}$ for $n=1, 2, \cdots$.  For $f(x)\in L^2(0, \pi)$, assume that there exists a limited number set $\Sigma_N=\{n_1, n_2, \cdots, n_N\}\subset {\bf N}$ such that
$$
(f(x), \sin(n x))\left\{\begin{array}{lll}
\neq 0, n\in \Sigma_N,\\
=0, \hbox{others},
\end{array}
\right.
\eqno{(4.1)}
$$
which means that there holds
$$
f(x)=\sum\limits_{n\in \Sigma_N} f_n\sin(nx),\eqno{(4.2)}
$$
where $f_n=\frac{2}{\pi} (f(x), \sin(nx))$.\\
\indent By (3.5) we have
$$
u(x, t)=\sum\limits_{n\in \Sigma_N} f_n E_{\alpha}(-D \lambda_n t^\alpha)\sin(nx),\ 0<t\leq T.
\eqno{(4.3)}
$$
Then the equation (3.7) reduces to
$$
F(\alpha)=\sum\limits_{n\in \Sigma_N} f_n E_{\alpha}(-D \lambda_n t_1^\alpha)\ \varphi_n(x_0)=d,
\eqno{(4.4)}
$$
and we get by using (2.8) in Lemma 3
$$
F'(\alpha)=\sum\limits_{n\in \Sigma_N} {f_n}\varphi_n(x_0) \sum\limits_{j=1}^\infty
\frac{(-1)^j\ j\ D^j \lambda_n^j t_1^{\alpha j}}{\Gamma(\alpha j+1)} [\ln(t_1)+\gamma_j],\ 0<\alpha<1.
\eqno{(4.5)}
$$
\indent Henceforth, by using Corollary 1 there holds
\begin{thm}
For suitable $x_0\in \Omega$ and $t_1>0$, there exists $\alpha_0\in I$ such that $F'(\alpha_0)\neq 0$, and then the equation (4.4) has a unique solution in a neighborhood of $\alpha_0$, and the solution is continuously dependent upon the measured data locally, i.e., the inverse fractional order problem is of local well-posedness.
\end{thm}
{\emph{Proof}}\quad By suitably choosing the measured point $x_0$, we have $\varphi_n(x_0)\neq 0$ for $n\in \Sigma_N$, and then there exists one $\alpha_0\in I$ such that $F'(\alpha_0) \neq 0$.
Then by using Corollary 1, the nonlinear equation (4.4) has a unique solution at the neighborhood of $\alpha_0$, i.e., the inverse fractional order problem is locally solvable. Moreover, noting the local continuity of the inverse function $F^{-1}$, the inverse problem is of local stability.\\
\indent Furthermore, we can get a conditional uniqueness under the nonnegative conditions for the initial values.\\
\subsection{Conditional uniqueness}
By Theorem 1 the nonlinear equation (4.4) is locally solvable for $d=u(x_0, t_1)$  where $0<t_1\leq T$ and $x_0\in \Omega$ be fixed. There are no global uniqueness for solving nonlinear equations in general case, but the situation could be changed if special conditions are overposed by which a conditional uniqueness is obtained. For convenience we also set $\Omega=(0, \pi)$. There holds
\begin{thm}
Assume that the the initial value function has the form of (4.2), and there hold $f_n \varphi_n(x)\geq 0$ for $n\in \Sigma_N$ and $x\in \Omega$. Then the fractional order can be uniquely determined by one measurement $u(x_0, t_1)$ if $\varphi_{n}(x_0)\neq 0$ and $t_1>0$.
\end{thm}
{\emph{Proof}}\quad Under the given conditions, there are $\varphi_n(x)=\sin(n x)$ and $\lambda_n=n^2$ for $n=1, 2, \cdots$.
By (4.3) the solution of the forward problem is given as
$$
u(\alpha)(x, t)= \sum\limits_{n\in \Sigma_N} f_n E_\alpha(-D n^2 t^\alpha) \sin(n x),
\eqno{(4.6)}
$$
and we get the nonlinear equation on the order $\alpha\in (0, 1)$
$$
F(\alpha)=\sum\limits_{n\in \Sigma_N} E_\alpha(-D n^2 t_1^\alpha)\  f_n \sin(n x_0)=d,\eqno{(4.7)}
$$
where $d=u(x_0, t_1)$ is the additional measurement.\\
\indent By choosing $x_0\in \Omega$ such that $\sin(n x_0)>0$ for $n\in \Sigma_N$, and noting the condition $f_n \varphi_n(x)\geq 0$ for $n\in \Sigma_N$, we can deduce that the function $F(\alpha)$ is of strict monotonicity on $\alpha\in (0, 1)$ due to the strict monotonicity and positivity of the Mittag-Lelffer function. Therefore the equation (4.7) has only one solution $\alpha^*\in I=(0, 1)$. The proof is completed.
\subsection{Theoretical testification}
\subsubsection{Example 1}
In this example we set $D=0.1$, and the initial function $f(x)=\frac{1}{2} \sin(2 x)$ for $x\in [0, \pi]$, and we choose $x_0=\frac{\pi}{4}$ as the measured point, and $t_1=2$ as the measured time, i.e., the additional measurement is $d=u(\frac{\pi}{4}, 2)$. Noting $n=2$ and $f_n=\frac{1}{2}$, and $\sin(n x_0)=1$, we have by (4.7)
$$
F(\alpha)=\frac{1}{2} {E_\alpha}(-\frac{2^{1+\alpha}}{5}),\ 0<\alpha<1,
\eqno{(4.8)}
$$
and we need to solve the equation
$$
F(\alpha)-d = 0,\eqno{(4.9)}
$$
using the additional data $u(x_0, t_1)=d$.\\
\indent Let the exact order be $\alpha^{exa}=0.75$, by which we get the additional data $d=0.25818$ by solving the forward problem. In order to see the unique solvability of the equation, we plot the continuous function $F(\alpha)-d$ on $\alpha\in [0, 1]$ in Figure 1, where $F(0)$ and $F(1)$ are defined by
$$
F(0)=\frac{1}{2}\sum\limits_{k=0}^\infty (-\frac{2}{5})^k,\ F(1)=\frac{1}{2} e^{-\frac{4}{5}}.\eqno{(4.10)}
$$
It can be seen clearly that the function $F(\alpha)-d$ is strictly monotone on $\alpha\in [0, 1]$, and the inverse order problem is of uniqueness.

\vbox{
\begin{center}
\includegraphics[scale=0.85]{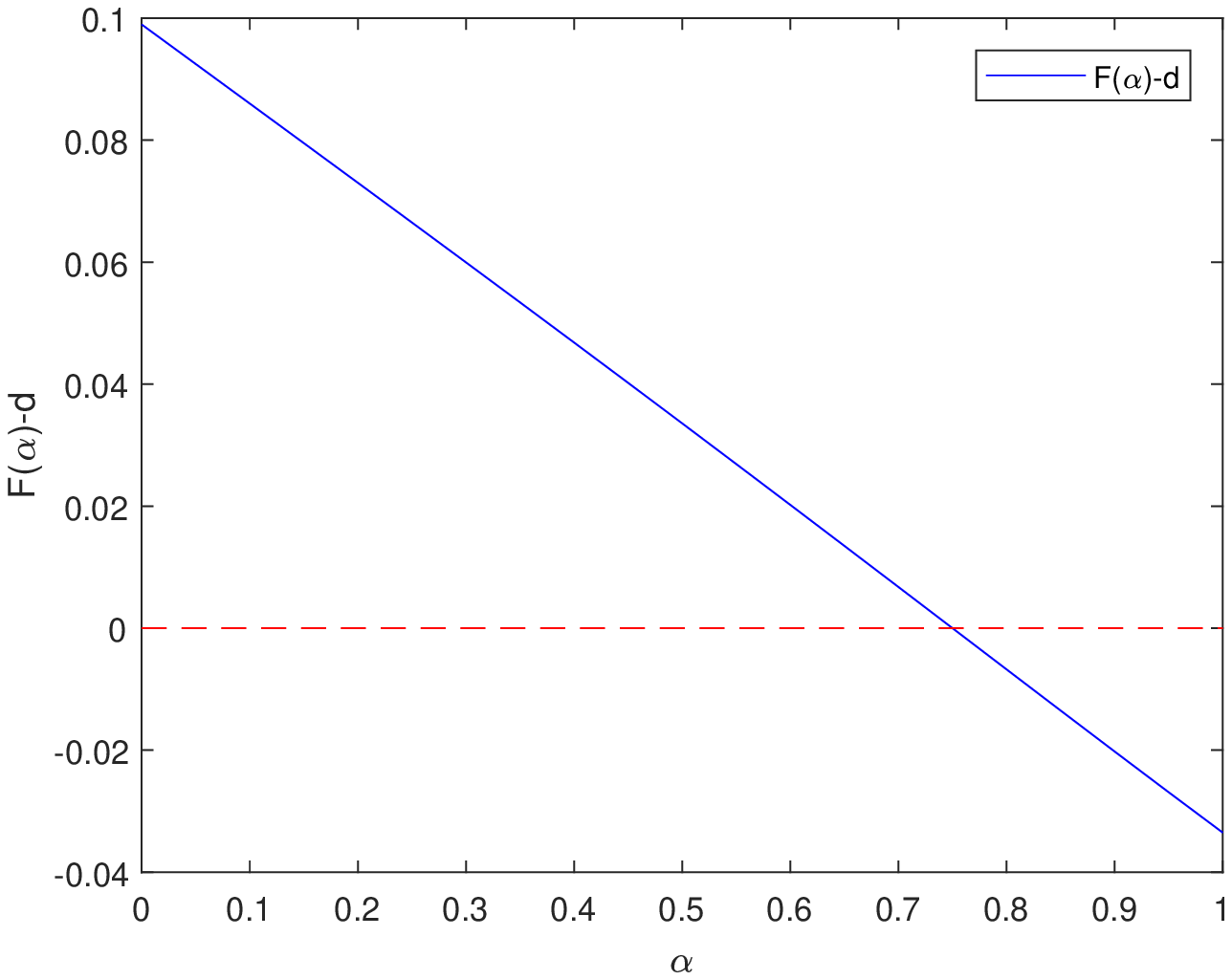}\\
\small{Figure 1. The picture of $F(\alpha)-d$ on $\alpha\in [0, 1]$ in Ex. 1.}
\end{center}
}

\subsubsection{Example 2}
We set $D=0.05$, and the initial function $f(x)=2 \sin(x)+\frac{1}{2} \sin(3 x)$ for $x\in [0, \pi]$ in this example. For given $\alpha\in (0, 1)$, noting that $D=0.05$, and $\lambda_1=1, \lambda_3=3$, the solution of the forward problem in this example is expressed as follows
$$
u(\alpha)(x, t)= 2 E_\alpha(-0.05 t^\alpha) \sin(x)+\frac{1}{2} E_\alpha(-0.45 t^\alpha) \sin(3x).
\eqno{(4.11)}
$$
\indent For the inverse order problem, let the exact fractional order be $\alpha^{exa}=0.5$. We choose $x_0=\frac{\pi}{6}$ as the measured point, and $t_1=10$ as the measured time, i.e., the additional measurement is $d=u(\frac{\pi}{6}, 10)=1.0112$. We have the nonlinear equation combing with the additional data
$$
F(\alpha)=E_\alpha(-0.05\cdot 10^\alpha)+\frac{1}{2} E_\alpha(-0.45\cdot 10^\alpha)=d.
\eqno{(4.12)}
$$
As done in the above, the function $F(\alpha)-d$ on $\alpha\in [0, 1]$ are plotted in Figure 2.

\vbox{
\begin{center}
\includegraphics[scale=0.85]{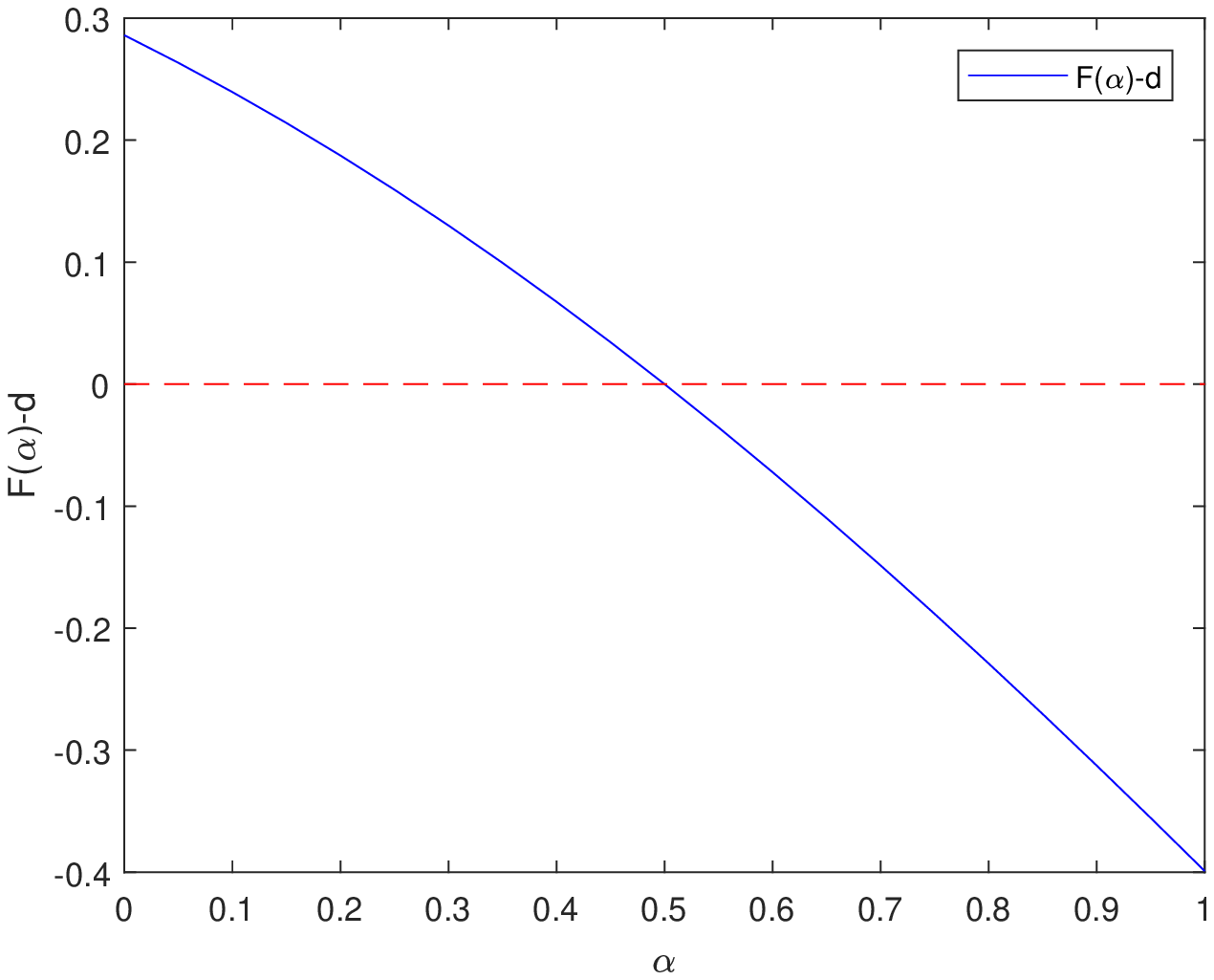}\\
\small{Figure 2. The picture of $F(\alpha)-d$ on $\alpha\in [0, 1]$ in Ex. 2.}
\end{center}
}

\indent From Figure 2 it can be seen again that the equation $F(\alpha)-d=0$ has a unique solution, and the inverse order problem is uniquely solvable by one measurement.

\section{Conclusion}
An inverse problem of identifying the fractional order in the 1D time fractional diffusion equation by using one measurement is investigated. Based on the expression of the solution to the forward problem, the inverse problem is reduced to a nonlinear algebraic equation on the fractional order, and the unique solvability can be obtained by the monotonicity and positivity of the Mittag-Lellfer function of real variables under the positive conditions for the initial values and the measured point. Theoretical examples are given to illustrate the conditional uniqueness of the inverse problem. It is noted that the derivative of the function $F(\alpha)$ on $\alpha\in (0, 1)$ can be computed by (4.8) or (4.12) respectively, and some gradient-type iterative algorithms can be applied to solve the nonlinear equation for which we will give details in another work in the near future.

\section*{Acknowledgements}
This work is supported by Natural Science Foundation of Shandong Province, China (No. ZR2019MA021), and National Natural Science Foundation of China (nos. 11871313, 11371231).



\begin{thebibliography}{99}

\bibitem {ada92} E. E. Adams and L. W. Gelhar,
Field study of dispersion in a heterogeneous aquifer 2: Spatial moments analysis,
Water Resources Research, 28(1992) 3293--3307.

\bibitem {ali15} S. Alimov and R. Ashurov,
Inverse problem of determining an order of the Riemann-Liouville time-fractional derivative,
Prog. Fractional Differ. Appl., 1(2015) 1--7.

\bibitem {ali20} S. Alimov and R. Ashurov,
Inverse problem of determining an order of the Caputo time fractional derivative for a subdiffusion equation,
Journal of Inverse and Ill-Posed Prolems, 28(2020) 651--658.

\bibitem {ash201} R. Ashurov and S. Umarov,
Determination of the order of fractional derivative for subdiffusion equations,
Fractional Calculus and Applied Analysis, 23(2020) 1647--1662.

\bibitem {ash202} R. Ashurov and R. Zunnunov,
Initial-boundary value and inverse problems for subdiffusion equations in ${\bf R}^N$,
Fractional Differ. Calculus, 10(2020) 291--306.

\bibitem {ben98} D. A. Benson,
The Fractional Advection-Dispersion Equation: Development and Application,
Dissertation of Doctorial Degree, University of Nevada, Reno, USA, 1998.

\bibitem {ber00} B. Berkowitz, H. Scher and S. E. Silliman,
Anomalous transport in laboratory-scale heterogeneous porus media,
Water Resources Research, 36(2000) 149--158. 

\bibitem {chen12} Z. Q. Chen, M. M. Meerschaert and E. Nane,
Space-time fractional diffusion on bounded domains,
J. Math. Anal. Appl., 393(2012) 479--488.

\bibitem {cheng09} J. Cheng, J. Nakagawa, M. Yamamoto and T. Yamazaki,
Uniqueness in an inverse problem for a one-dimensional fractional diffusion equation,
Inverse Problems, 25(2009) 115002. 

\bibitem {chi11} G. S. Chi, G. S. Li and X. Z. Jia,
Numerical inversions of source term in FADE with Dirichlet boundary condition by final observations,
Computers and Mathematics with Applications, 62(2011) 1619-1626.

\bibitem {ding15}
X. L. Ding and J. J. Nieto,
Analytical solutions for the multi-term time-space fractional reaction-diffusion on equations on an infinite domain,
Fractional Calculas and Applied Analysis, 18(2015) 697--716.

\bibitem {eide04} S. D. Eidelman and A. N. Kochubei,
Cauchy problem for fractional diffusion equations,
J. Diff. Eqn., 199(2004) 211--255.

\bibitem {gor99} R. Gorenflo, Y. Luchko and P. P. Zabrejko,
On solvability of linear fractional differential equations in Banach spaces,
Fract. Calc. Appl. Anal., 2(1999)163--176.

\bibitem {gor20} R. Gorenflo, A. A. Kilbas, F. Mainardi, and S. Rogosin,
Mittag-Leffler Functions, Related Topics and Applications (2nd Ed.),
Springer, Berlin, 2020.

\bibitem {gue21} N. Guerngar, E. Nane, R. Tinaztepe, S. Ulusoy and H. W. Van Wyk,
Simultaneous inversion for the fractional exponents in the space-time fractional diffusion equation,
Fractional Calculus and Applied Analysis, 24(2021).

\bibitem {han02} A. Hanyga,
Multidimensional solutions of time-fractional diffusion-wave equations,
Proc. R. Soc. Lond. A, 458(2002)933--957.

\bibitem {hat98} Y. Hatano and N. Hatano,
Dispersive transport of ions in column experiments: an explanation of long-tailed profiles,
Water Resources Research, 34(1998) 1027--1033.

\bibitem {hat13} Y. Hatano, J. Nakagawa, S. Wang, and M. Yamamoto,
Determination of order in fractional diffusion equation,
J. Math. Ind., 5A(2013) 51--57.

\bibitem {jan16} J. Janno,
Determination of the order of fractional derivative and a kernel in an inverse problem for a genaralized time fractional diffusion equation,
Electron. J. Differ. Equ., 199(2016) 1--28. 

\bibitem {jan18} J. Janno and N. Kinash,
Reconstruction of an order of derivative and a source term in a fractional diffusion equation from final measurements,
Inverse Problems, 34(2018) 025007. 

\bibitem {jia14} X. Z. Jia, D. L. Zhang, G. S. Li, et al,
Numerical inversion of the fractional orders in the space-time fractional advection-dispersion equation with variable coefficients (in Chinese).
Mathematica Numerica Sinica, 36(2014) 113--132.

\bibitem {jin15} B. T. Jin and W. Rundell,
A tutorial on inverse problems for anomalous diffusion processes,
Inverse Problems, 31(2015) 035003. 

\bibitem {jin21} B. T. Jin and Y. Kian,
Recovery of the order of derivation for fractional diffusion equations in an unknown medium,
2021, arXiv: 2101.09165.

\bibitem {kil06} A. A. Kilbas, H. M. Srivastava and J. J. Trujillo,
Theory and Applications of Fractional Differential Equations, Elsevier, Amsterdam, 2006.

\bibitem {kub20} A. Kubica, K. Ryszewska,  M. Yamamoto,
Theory of Time-Fractional Differential Equations an Introduction,
Springer, Berlin, 2020.

\bibitem {li13} G. S. Li, D. L. Zhang, X. Z. Jia and M. Yamamoto,
Simultaneous inversion for the space-dependent diffusion coefficient and the fractional order in the time-fractional diffusion equation,
Inverse Problems, 29(2013) 065014. 

\bibitem {liq05} Q. Y. Li, Z. Z. Mo, L. Q. Qi,
Numerical Solutions for Nonlinear Equations (in Chinese),
Science Press, Beijing, 2005.

\bibitem {liz151}
Z. Y. Li and M. Yamamoto,
Uniqueness for inverse problems of determining orders of multi-term time-fractional derivatives of diffusion equation,
Applicable Analysis, 94(2015) 570--579.

\bibitem {liz16}
Z. Y. Li, O. Y. Imanuvilov and M. Yamamoto,
Uniqueness in inverse boundary value problems for fractional diffusion equations,
Inverse Problems, 32(2016) 015004.

\bibitem {liz20} Z. Y. Li, K. Fujishiro, G. S. Li,
Uniqueness in the inversion of distributed orders in ultraslow diffusion equations,
Journal of Computational and Applied Mathematics, 369(2020) 112564.

\bibitem {liuj10} J. J. Liu and M. Yamamoto,
A backward problem for the time-fractional diffusion equation,
Applicable Analysis, 89(2000) 1769--1788.

\bibitem {liuy16}
Y. K. Liu, W. Rundell and M. Yamamoto,
Strong maximum principle for fractional diffusion equations and an application to an inverse source problem,
Fractional Calculus and Applied Analysis, 19(2016) 888--906.

\bibitem {luch09} Y. Luchko,
Maximum principle for the generalized time-fractional diffusion equation,
J. Math. Anal. Appl., 351(2009)218--223.

\bibitem {luch10} Y. Luchko,
Some uniqueness and existence results for the initial-boundary-value problems for the generalized time-fractional diffusion equation,
Comput. Math. Appl. 59(2010)1766--1772.

\bibitem {main10} F. Mainardi,
Fractional Calculus and Waves in Linear Viscoelasticity: An Introduction to Mathematical Models,
Imperial College Press, London, 2010.

\bibitem {metz00} R. Metzler and J. Klafter,
Boundary value problems for fractional diffusion equations,
Physica A, 278(2000) 107--125.

\bibitem {mur07} D. A. Murio,
Stable numerical solution of fractional-diffusion inverse heat conduction probllem,
Computers and Mathematics with Applications, 53(2007) 1492--1501.

\bibitem {pod99} I. Podlubny,
Fractional Differential Equations, Academic, San Diego, 1999.

\bibitem {saka11} K. Sakamoto and M. Yamamoto,
Initial value/boundary value problems for fractional diffusion-wave equations and applications to some inverse problems,
Journal of Mathematical Analysis and Applications, 382(2011) 426--447.

\bibitem {sunl21}
L. L. Sun, Y. S. Li and Y. Zhang,
Simultaneous inversion for the potential term and the fractional orders in a multi-term time-fractional diffusion equation,
Inverse problems, 37(2021) 055007.

\bibitem {tatar13} S. Tatar and S. Ulusoy,
A uniqueness result for an inverse problem in a space-time fractional diffusion equation,
Electronic Journal of Differential Equations, 258(2013): 1--9.

\bibitem {weit16}
T. Wei, X. L. Li and Y. S. Li,
An inverse time-dependent source problem for a time-fractional diffusion equation,
Inverse problems, 32(2016) 085003.

\bibitem {xiong06} Y. W. Xiong, G. H. Huang and Q. Z. Huang,
Modeling solute transport in one dimensional homogeneous and heterogeneous soil columns with continuous time random walk,
J. Contam. Hydrol., 86(2006) 163-175. 

\bibitem {yama12} M. Yamamoto,  Y. Zhang,
Conditional stability in determining a zeroth-order coefficient in a half-order fractional diffusion equation by a Carleman estimate,
Inverse Problmes, 28(2012) 105010.

\bibitem {yama20} M. Yamamoto,
Uniqueness in determining the orders of time and spatial fractional derivatives,
2020, arXiv: 2006.15046.

\bibitem {yama21} M. Yamamoto,
Uniqueness in determining fractional orders of derivatives and initial values,
Inverse Problems, 37(2021) 095006.

\bibitem {yu15} B. Yu, X. Jiang and H. Qi,
An inverse problem to estimate an unknown order of a Riemann-Liouville fractional derivative for a fractional Stokes' first problem
for a heated generalized second grade fluid,
Acta Mech. Sin., 31(2015) 153--161.

\bibitem {zhou03} L. Zhou and H. M. Selim,
Application of the fractional advection-dispersion equations in porous media,
Soil. Sci. Soc. Am. J., 67(2003) 1079--1084.

\end {thebibliography}
\end{document}